\documentclass[preprint,12pt]{elsarticle}

\usepackage{graphicx}
\usepackage{psfrag} 
\usepackage{mathptmx}      
\usepackage{amssymb}
\usepackage{amsmath}
\usepackage{amsthm}
\usepackage[dvipsnames]{xcolor}
\usepackage{natbib}
\usepackage{marvosym}
\usepackage{a4wide}

\usepackage[mathscr]{euscript}

\newtheorem{theorem}{Theorem}[section]
\newtheorem{lemma}[theorem]{Lemma} 

\newtheorem{proposition}[theorem]{Proposition}

\newtheorem{observation}[theorem]{Observation}

\newtheorem{remark}[theorem]{Remark}

\newcommand{\pf}{\noindent{\em Proof: }}
\newcommand{\epf}{\hfill\hbox{\rule{3pt}{6pt}}\\}

\journal{Advances in Applied Mathematics}

\begin{document}
	
	\begin{frontmatter}

\author{K.T. Huber
	\fnref{label2}}
\ead{k.huber@uea.ac.uk}
\author{V. Moulton\fnref{label2}}
\ead{v.moulton@uea.ac.uk}
\author{G. E. Scholz\fnref{label3}}

\cortext[cor1]{Corresponding author: K.T. Huber}
\fntext[label2]{School of Computing Sciences, 
University of East Anglia, UK}
\ead{gllm.scholz@gmail.com}
\address{\fnref{label2}}
\fntext[label3]{Laboratoire d'Informatique, de Robotique et de
	Micro\'el\'ectronique de
	Montpellier (LIRMM), Universit\'e de Montpellier,
	Montpellier, France.}
\address{\fnref{label3}}

\title{Overlaid species forests}




\begin{abstract}
Introgression  is an evolutionary process in which 
genes or other types of genetic material are introduced into a genome. 
It is an important evolutionary process that can, for example, play a fundamental
role in speciation. Recently the concept of an overlaid species
forest was introduced to represent introgression histories. Basically
this approach takes a putative gene history in the form
of a phylogenetic gene tree and tries to overlay this onto a forest which usually
consists of a collection of lineage trees for the species of interest.
The result is a network called an overlaid species forest
in which genes jump or introgress between lineages.
In this paper we study properties of overlaid species
forests, showing that they have various connections
with models for lateral gene transfer, maximum parsimony, and unfolding of
phylogenetic networks. In particular, 
we show that a certain algorithm called \textsc{OSF-Builder} 
for constructing overlaid species forests is 
guaranteed to a produce a special type of overlaid species forest 
with a minimum number introgressions, as
well as providing some characterizations for networks that
can arise from  overlaid species forests. We expect
that these results will be useful in developing new methods
for representing introgression histories, a growing area 
of interest in phylogenetics. 
\end{abstract}


\end{frontmatter}

\section{Introduction}

    Introgression is an evolutionary process in which 
	genes or other types of genetic material are introduced into a genome \cite{goulet2017hybridization}.
	It is an important process since it can, for example, help species adapt to or 
	expand into new environments \cite{grant2016introgressive}.  Introgression
	is a wide-spread phenomenon in plants and animals, and 
	it often arises when species belonging to distinct lineages 
	of the same species hybridize \cite{mallet2005hybridization}. Some striking examples of
	introgression include  genes
	introduced by Neanderthals into modern humans  \cite{sankararaman2014genomic} and 
	genetically modified crop genes moving into their wild relatives \cite{stewart2003transgene}.
	Other examples are given in e.g. \cite{mallet2005hybridization,zhang2016genome}.
	Although several approaches have been introduced to 
	detect introgression (e.g. the ABBA/BABA
	test  \cite{sousa2013understanding}), to date relatively few approaches have
	been proposed for reconstructing explicit evolutionary scenarios 
	which involve introgression (see \cite{solis2016inferring,wen2016reticulate} for some examples).
	 
\begin{figure}[h]
		\begin{center}
			\includegraphics[scale=0.5]{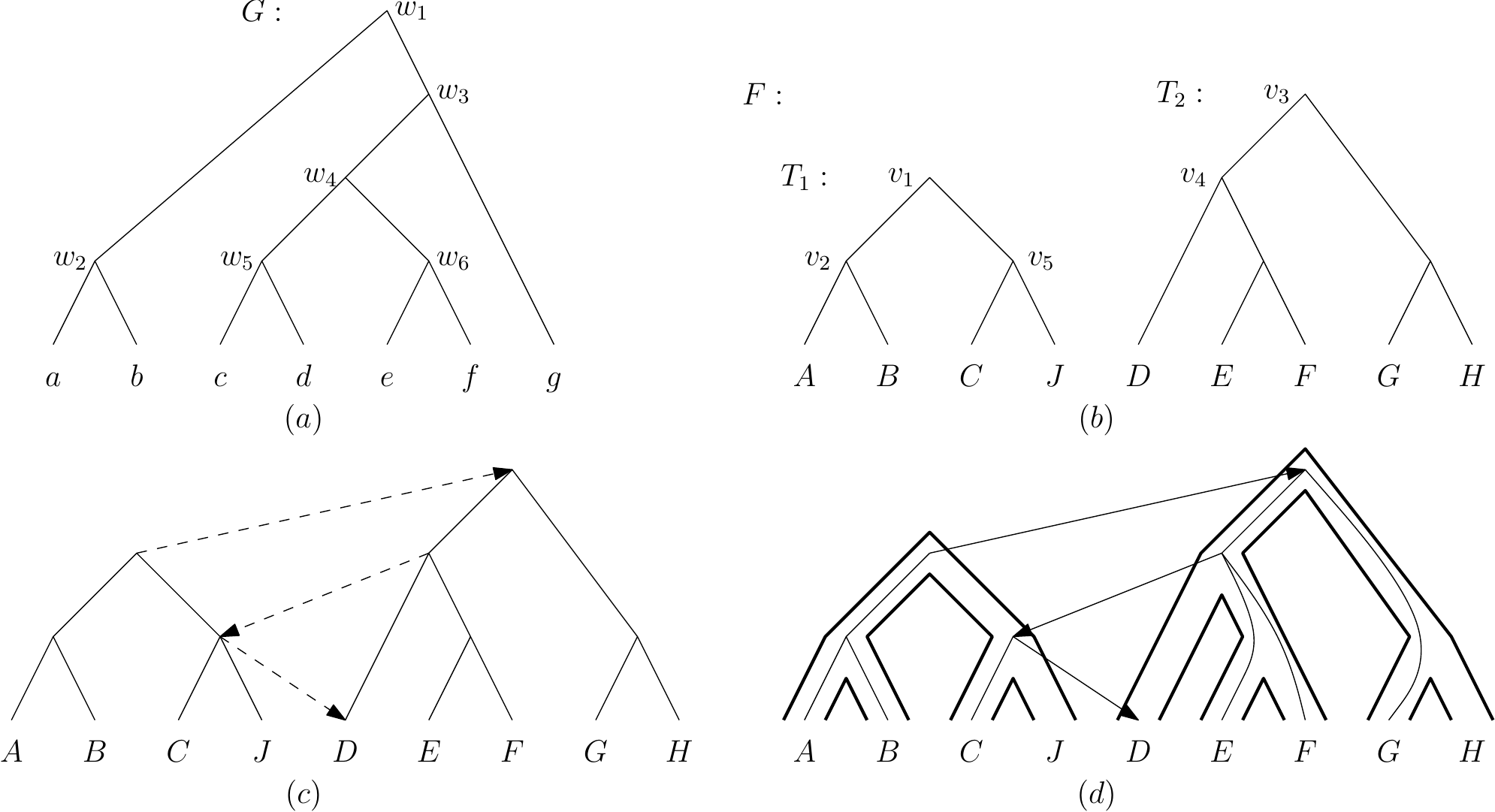}
			\caption{A forest triple where $G$ is the
                            tree in (a), $F$ is the
                            forest $\{T_1,T_2\}$ in (b) and $\phi$ is the
                            map that takes every leaf of $G$ to the
                            corresponding leaf of $F$ in capitals. (c) A representation of an OSF $\psi$ for $(G,F,\phi)$
                            where 
                            $\psi(w_i)=v_i$ for $i=1,2,3,5$ and $\psi(w_6)=\psi(w_4)=v_4$.
                            The dashed edges are contact arcs.
                            This representation is visualised in more detail in (d).
                          \label{fig:introduce}}
		\end{center}
	\end{figure}

    Recently, a new     
    method was introduced for constructing introgression histories
    which involves overlaying a species forest $F$ with a gene tree $G$ \cite{scholz2019osf}.
    In this method, $F$ usually represents a collection of lineage trees for a 
	certain species (e.g. butterflies), and $G$ the evolutionary 
	history of some gene which 
	jumps or introgresses between the 
	lineages (e.g. genes which affect butterfly wing colouring). 
	Both are linked via the map $\phi$ from the
	leaf-set of $G$ to the leaf-set of $F$ which is
	given by taking a gene to the species  in which it is 
	resides, giving 
	a triple $(G,F,\phi)$ which we call a a {\em forest triple}.
	For example, in Figure~\ref{fig:introduce} we present 
	a gene tree $G$ and a forest $F$ consisting of 
	two trees $T_1$ and $T_2$
	together with the map $\phi$ which takes lower case labels on $G$
	to their labels in capitals in $F$ (e.g. $\phi(b)=B$). 
	Based on this information, we  
	then look for ways to overlay $G$ onto the forest $F$, i.e. map
	the vertex set of $G$ to the vertex set of $F$ whilst respecting $\phi$ and other conditions of ancestry,  
	so as to represent the introgression history of the gene in question.
	For example, in Figure~\ref{fig:introduce}(c), we present one possible
	mapping $\psi$ of the vertex set of the depicted gene tree
	$G$ to the vertex set of the forest $F$ also depicted in the figure.
	We call such a mapping an {\em overlaid species forests} or OSF for short.

	In this paper, we shall study properties
	of OSFs. We begin by introducing some 
	notation in Section~\ref{sec:prelims}, and reviewing the \textsc{OSF-Builder} algorithm.
	This algorithm is introduced in \cite{scholz2019osf} and, for a given forest triple, aims to produce 
	an OSF that minimizes the number of times that a gene introgresses between lineages.
	The rationale behind this is that introgression is thought to be relatively infrequent in nature.
    In Section~\ref{sec:strictosf}, we show that the 
	\textsc{OSF-Builder} algorithm is guaranteed to produce an 
	optimal OSF, and that this OSF is of a special form, called a \emph{strict OSF} (Theorem~\ref{th:minimum}).
	This type of OSF is related to a well-studied tree reconciliation model for representing lateral gene transfer for trees
	introduced in \cite{tofigh2011simultaneous}, the main 
	difference being that in our model 
	a tree is replaced by a forest. Our approach to proving 
	this result uses the close relationship between the \textsc{OSF-Builder} algorithm and
	the Fitch-Hartigan algorithm for computing most parsimonious trees \cite{F71,H73}. 
	
	In Section~\ref{sec:valid} we then look more deeply into structural properties of OSFs.
	Since each OSF can be visualized in terms of a network, as depicted in e.g. Figure~\ref{fig:introduce}(c), 
	it is thus of interest to find characterizations for those networks 
	which arise from OSFs. We call such networks \emph{valid} networks. 	Note that this is related to the problem of
	characterizing tree-based networks \cite{francis2015phylogenetic}, a problem
	that has recently received a great deal of attention in the
	literature (see Section~\ref{sec:conclusion} for more details). In Section~\ref{sec:valid}
	we show that a network is a valid network if and only 
	if it can be obtained by applying 
	\textsc{OSF-Builder} algorithm to some forest triple (Theorem~\ref{thm:valid}).
	Then, in Section~\ref{sec:representing} we present a 
	structural characterization of valid networks (Theorem~\ref{thm:structure}). To prove this result, we 
	use a variation of the process of ``unfolding" a phylogenetic network \cite{huber2016folding}.

    In Section~\ref{spr-forest}, we consider the affect
	of changing the gene tree and species forest has
	on the number of arcs contained in an optimal network 
	corresponding to a forest triple. This is of interest because
	it helps to shed light on how noise in a 
	forest triple can affect the OSF's generated by the  \textsc{OSF-Builder} algorithm.
	In particular, we provide upper bounds on the amount that the optimal 
	score for an OSF can change by
	in terms of the number of  \emph{subtree prune and regraft (SPR) operations} 
	that alter a gene tree or forest  (Theorems~\ref{spr-tree} and \ref{spr-forest}). 
	We do this by exploiting some results in \cite{FK16}
	which relate the so-called maximum parsimony and SPR-distances 
	between two phylogenetic trees.
	In the last section we conclude by briefly discussing  some possible future directions.

\section{Preliminaries}
\label{sec:prelims}

In what follows, we assume $X$ is a finite set with $|X|\ge 2$ 
and that all graphs are directed unless otherwise stated, and without loops.

\subsection{Networks} 

Let $G$ be a directed graph. We denote the vertex
set of $G$ by $V(G)$ and its set of arcs by $A(G)$.
We denote an arc $a$ from a vertex $u\in V(G)$
to a vertex $v\in V(G)$ by $a=(u,v)$ and refer to
$u$ as $tail(a)$ and to $v$ as $head(a)$, respectively.
A {\em walk} in $G$ is a sequence of vertices $v_1,\dots,v_m$ in $V(G)$,
$m \ge 1$, 
such that $(v_i,v_{i+1}) \in A(G)$, $1 \le i \le m-1$.
A {\em trail} in $G$ is a walk with no repeated arcs, and a
{\em path} is a trail with no repeated vertices.
A {\em cycle} in $G$ is a trail of the form $v_1,\dots,v_m,v_1$, $m \ge 2$
where $v_1,\dots,v_m$ is a path.
We say that $G$ is {\em connected} if its underlying (undirected)
graph is connected (note that the underlying graph might contain
multi-edges).  

We denote by $outdeg_G(v)=outdeg(v)$
the number of outgoing arcs of a vertex $v\in V(G)$
and by $indeg_G(v)=indeg(v)$ its number
of incoming arcs. 
We call a vertex $v\in V(G)$ with $indeg(v)=1$ and $outdeg(v)=0$
a \emph{leaf} of $G$ and denote by $L(G)$ the set of leaves of $G$.
We call a vertex $v\in V(G)$ with $indeg(v)=0$ and $outdeg(v)\ge 2$ 
a \emph{root} of $G$.

A {\em network (on $X$)} is a directed graph $G$ which satisfies:
\begin{itemize}
\item[(N1)] $X \subseteq V(G)$,
\item[(N2)] if $v\in V(G)$ such that $indeg(v)=1=outdeg(v)$
  then $v\in X$,
\item[(N3)] $\{v \in V(G) \,: \, indeg(v)=0\} \cap X = \emptyset$, and
\item[(N4)] $\{v \in V(G) \,: \, outdeg(v)=0\}  \subseteq X$ 
\end{itemize}

Two networks $N$, $N'$ on $X$ are {\em isomorphic} if 
there exists a bijective map $\psi:V(N)\to V(N')$
that induces a graph isomorphism between $N$ and $N'$ that is the
identity on $X$.

Suppose $G$ is a network and
$u,v\in V(G)$. Then we put $u\preceq_G v$ (or just $u\preceq v$)
if there is a directed path in $G$ starting at $u$ and ending at $v$. 
If $u\preceq_G v$ then we call $u$ an \emph{ancestor} 
of $v$ (in $G$), and say that $v$ lies \emph{below} $u$ (in $G$). 
Note that a vertex can be its own ancestor. A vertex $u$ is a {\em child}
of a vertex $v$ if $(v,u) \in A(G)$. 

A {\em (rooted) phylogenetic tree $T$ (on $X$)}
is  a network $T$ whose underlying graph is a tree, that has a single root
denoted $\rho_T$,  and leaf set $X$ (and so $L(T)$ has size at least 2,
by our assumption on $X$).
We denote the set of leaves of $T$ below  $u$
by $\mathcal C_T(u)$. For any non-empty subset $Y\subseteq X$,
we denote by $lca_T(Y)$ the unique vertex $v \in V(T)$ that is an
ancestor of every element in $Y$  such that no vertex 
below $v$ and distinct from $v$ is an ancestor of every element of $Y$.
If $Y=\{y_1,\ldots, y_k\}$, $k\geq 1$, then we sometimes write
$ lca_T(y_1,\ldots, y_k)$ rather than $lca_T(\{y_1,\ldots, y_k\})$. Note that
$lca_T(x)=x$, for all $x\in X$.

\subsection{Overlaid species forests} \label{subsec:osf}

A {\em forest} $F$ is a non-empty set of phylogenetic trees.
Note that the leaf set of a forest $F$ is $\bigcup_{T\in F} L(T)$. 
To reduce notation, we shall also regard a forest as being a graph (whose
components are the trees in $F$),
and in case a forest has one element we shall also consider this
as being a phylogenetic tree.
A vertex in a forest $F$ that is not a leaf is called 
an {\em interior vertex} of $F$.  
We let $V^0(F)$ denote the set of interior vertices of $F$.
We say that a forest $F$ is {\em binary}
if $outdeg(v)=2$ for all interior vertices $v$ of $F$.

A triple $\mathcal F=(G,F,\phi=\phi_{G,F})$ consisting of
a phylogenetic tree $G$, a forest $F$, and a (leaf) map
$\phi=\phi_{G,F}:L(G)\to L(F)$ 
is called a {\em forest triple}. We say that $\mathcal F$ is
{\em binary} if both $G$ and $F$ are binary.
To ease readability of our examples, we usually denote the leaves of
$G$ by lower-case letters (also with indices), and the leaves in $F$ which 
they map to under $\phi$ by the corresponding capital letter.
Given a forest triple $\mathcal F=(G,F,\phi)$,
an {\em overlaid species forest or OSF (for $\mathcal F$)}, is  a map
$\psi: V(G)\to V(F)$ which satisfies:
\begin{itemize}
\item[(P1)] $\psi|_{L(G)}=\phi$.
\item[(P2)] If $u, v\in V(G)$ satisfy $u\preceq_G v$  and 
$\psi(u),\psi(v)\in V(T)$ holds for some $T\in F$ 
then $\psi(u)\preceq_T\psi(v)$.
\item[(P3)] If $u \in V(G)$ is such that $\psi(u) \in V(T)$ for some
 $T \in F$, then there exists some leaf $v \in \mathcal C_G(u)$ such that
 $\phi(v) \in \mathcal C_T(\psi(u))$.
\end{itemize}

For example, the map $\psi:V(G)\to V(F)$ given in the caption of
Figure~\ref{fig:introduce} is an OSF for 
the binary forest triple $\mathcal F=(G,F,\phi)$ given in that figure.

\subsection{The \textsc{OSF-Builder} algorithm}
\label{sec:osf-builder}

Suppose that $\psi$ is an OSF for a forest triple  $\mathcal F = (G,F,\phi)$.
A {\em contact arc of $\psi$} is a pair $(\psi(u),\psi(v))$
where $(u,v) \in A(G)$ 
and $\psi(u)\in V(T)$ and $\psi(v)\in V(T')$ for $T,T' \in F$ distinct. 
We let $C(\psi)$ denote the multi-set of contact arcs of $\psi$, 
with multi-set cardinality $|C(\psi)|$, and we set 
$$
t(\mathcal F) = min\{ |C(\psi)| \,:\, \psi \mbox{ an OSF for }  \mathcal F\}.
$$
In addition, we let $C^*(\psi)$ denote the underlying set of $C(\psi)$.
For example, for the OSF $\psi$ in
	Figure~\ref{fig:construction}, we have
	$|C(\psi)|=3$ and $|C^*(\psi)|=2$. 
In \cite{scholz2019osf} an algorithm called \textsc{OSF-Builder}
is introduced for computing an OSF $\psi$ for $\mathcal F$
with $|C(\psi)|=t(\mathcal F)$ which we now
briefly review. It is based on the Fitch-Hartigan 
algorithm for computing the parsimony score of a
character on a phylogenetic tree \cite{H73} (see also \cite{F71}).

Suppose $\mathcal F=(G,F,\phi)$ is a forest triple.
Consider the map $f = f_{\mathcal F}: L(G) \to F$, which takes
each $v \in L(G)$ to the tree $T \in F$ with $\phi(v) \in V(T)$, as
being a character\footnote{For an  arbitrary set $Y$, a {\em character on $Y$} 
is a map from $Y$ to some finite set.} on $L(G)$. 
An {\em extension $\bar{f}$ of $f$ to 
$V(G)$} is a map $\bar{f}: V(G) \to F$ such that
$\bar{f}(v)=f(v)$ for all $v \in L(G)$.
The {\em parsimony score $l_f(G)$ of 
$f$ on $G$} is then given by taking the minimum,
over all extensions $\bar{f}$ of $f$, of the number of
arcs $(u,v) \in V(G)$ with $\bar{f}(u) \neq \bar{f}(v)$ (cf. \cite[p.84]{semple2003phylogenetics}). 

\begin{figure}[h]
	\begin{center}
		\includegraphics[scale=0.7]{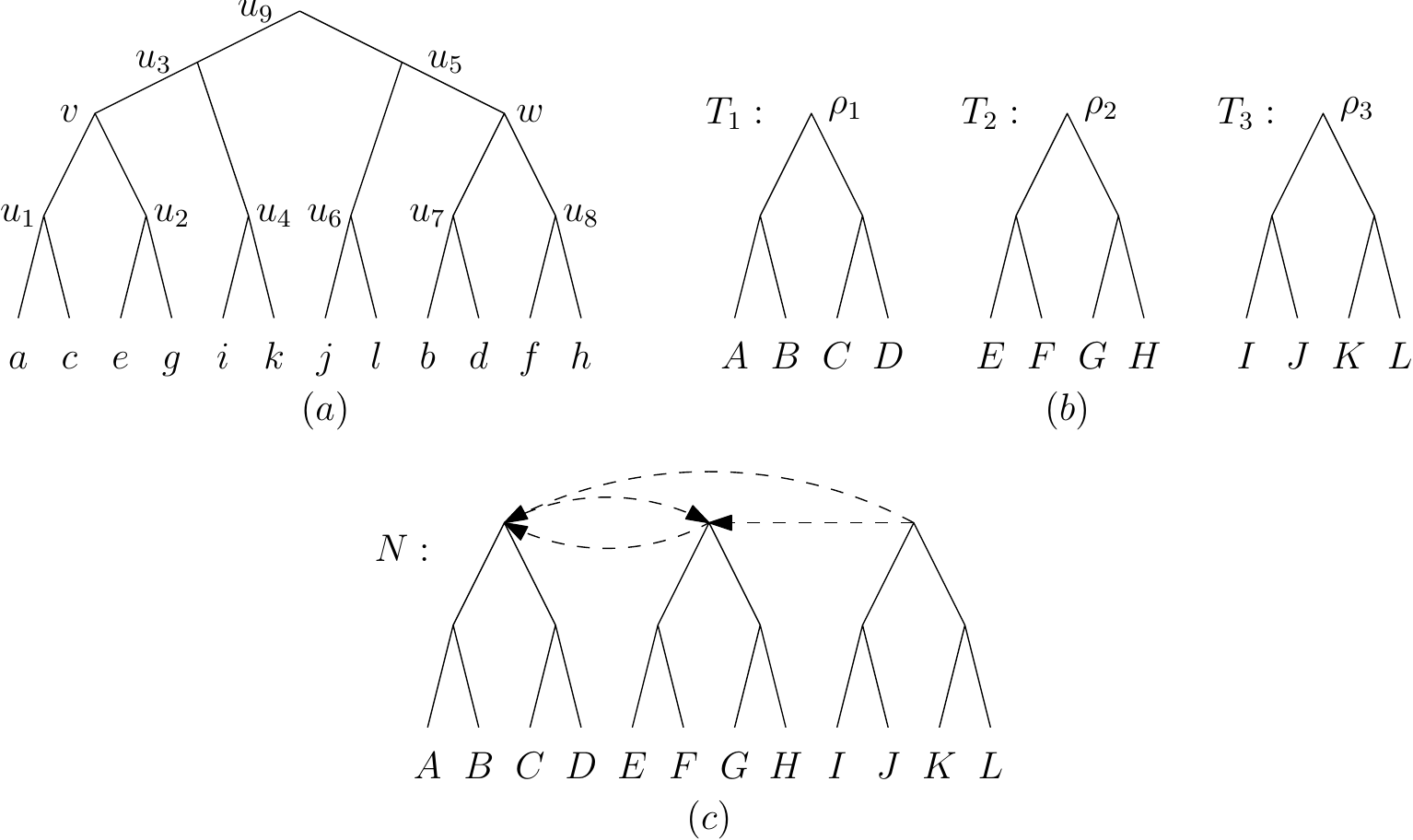}
		\caption{A forest triple $\mathcal F$ consisting of
			the tree $G$ depicted in (a), the forest 
			$F=\{T_1,T_2,T_3\}$ in (b), and leaf-map defined as usual. The network $N$ in (c) is the representation
			for the OSFs $\psi,\psi':V(G)\to V(F)$
			on $\mathcal F$ given by $\psi(u_i)=\psi'(u_i)=\rho_1$ if $i\in \{1,7\}$,
			$\psi(u_i)=\psi'(u_i)=\rho_2$ if $i\in \{2,8\}$,
			$\psi(u_i)=\psi'(u_i)=\rho_3$ otherwise, and
			$\psi(v)=\psi'(w)=\rho_2$ and $\psi(w)=\psi'(v)=\rho_1$.}
		\label{fig:partrep}
	\end{center}
\end{figure}

The \textsc{OSF-Builder} algorithm works by first 
associating to every vertex $v\in V^0(G)$
the set $\sigma(v)$ of trees in $F$
which are assigned most frequently to its children
in a bottom-up fashion, starting at the leaves.
Then it computes an OSF $\psi$ in a top-down phase as follows.
It begins by computing an extension $\bar{f}$ of $f$.
To do this it initializes first $\bar{f}(\rho_G)$ to
be any tree in $\sigma(\rho_G)$. 
Then, moving down $G$, for any vertex $v$ in $V^0(G)$ with $\bar{f}$
defined on $v$ but not yet on its children, for $u$ a child of $v$ 
it sets $\bar{f}(u) = \bar{f}(v)$ if 
$\bar{f}(v)$ is contained in $\sigma(u)$, and otherwise it 
sets $\bar{f}(u)$ to be any element in $\sigma(u)$.  The OSF $\psi$
is then defined by setting 
$$
\psi(v) = lca_{\bar{f}(v)}( \{w\in L(\bar{f}(v)) \,:\,
\mbox{ there exists some leaf }
x\in L(G) \mbox{ below } v \mbox{ such that } w=\phi(x) \}),
$$
for all $v \in V(G)$. In other words, $\psi(v)$ is  the last common ancestor of
all leaves in the tree $\bar{f}(v)$ that are the image
of some leaf in $\mathcal C_G(v)$ under $\phi$.
It follows that \textsc{OSF-Builder} computes an OSF $\psi$ on $\mathcal F$ 
such that $|C(\psi)|=l_{f_{\mathcal F}}(G)$. 
See Figure~\ref{fig:partrep} for an example of an OSF computed using  \textsc{OSF-Builder}.

\section{Strict OSFs}
\label{sec:strictosf}

\begin{figure}[h]
	\begin{center}
		\includegraphics[scale=0.7]{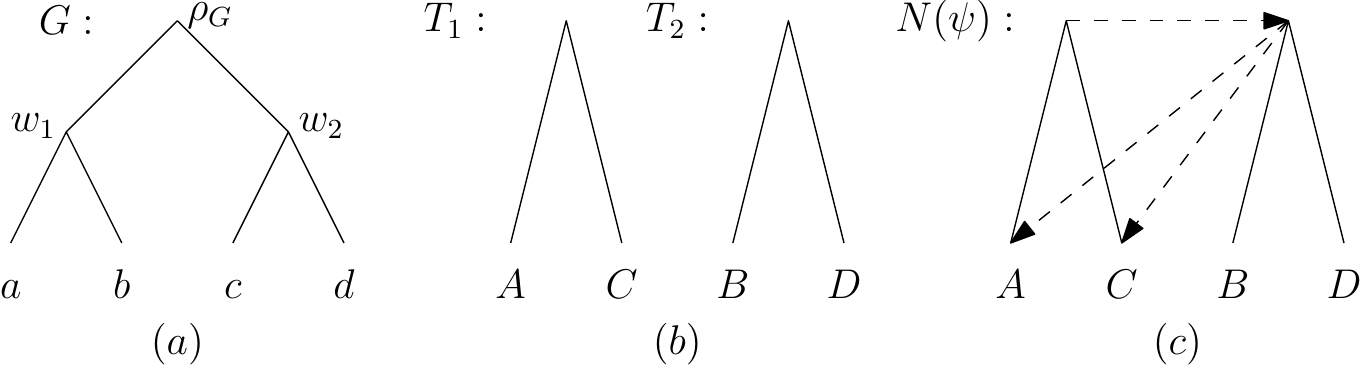}
		\caption{For $G$ the tree in (a),
			$F=\{T_1,T_2\}$ the forest in (b), and $\phi$ 
			the corresponding leaf-map, we present in (c)
			the representation $N(\psi)$
			of the OSF $\psi:V(G)\to V(F)$ for the forest triple $(G,F,\phi)$
			where $\psi$ maps the root of $G$ to the root of $T_1$ and $w_1$ and $w_2$ to
			the root of $T_2$. The dashed arcs are the contact arcs of $N(\psi)$. Note that $\psi$ 
			is not an sOSF for $\mathcal F$ because $\psi$
			does not satisfy (S3) for $\rho_G$.}
		\label{fig:osf-not-sosf}
	\end{center}
\end{figure}

In this section we show that the \textsc{OSF-Builder} algorithm produces 
an sOSFs  $\psi$ with $|C(\psi)|= t(\mathcal F)$ (see also \cite[Supp. Mat., Theorem 2]{scholz2019osf} 
for a related, weaker result). 
We begin by defining this special type of OSF. 
Suppose that  $\mathcal F = (G,F,\phi)$ is a forest triple.
A {\em strict overlaid species forest or sOSF (for $\mathcal F$)} is  a map
$\psi: V(G)\to V(F)$ which satisfies (P1), (P2) and 
\begin{itemize}
\item [(S3)]  If $u \in V^0(G)$ is such that $\psi(u) \in V(T)$
  for some
	$T \in F$, then there exists a child $v$ of $u$ in $G$ such that
	$\psi(v)$ is below $\psi(u)$ in $T$.
\end{itemize}
sOSFs are closely related to DTL-scenarios as 
defined in \cite[p.519]{tofigh2011simultaneous},
the main difference being that a DTL-scenario is essentially
a map from a (phylogenetic) gene tree to a species tree,
rather than to a forest.
Note that every sOSF for $\mathcal F$
is  an OSF for $\mathcal F$ (see Lemma~\ref{lem:trivial}) but not conversely
(see e.g. Figure~\ref{fig:osf-not-sosf}).

We now present some basic properties of sOSF's.

\begin{lemma}\label{lem:trivial}
	Suppose $\mathcal F=(G,F,\phi)$ is forest triple.
	\begin{itemize}
	\item[(i)] If $\psi$ is an sOSF for $\mathcal F$,
          then it is an OSF for $\mathcal F$.
	\item[(ii)] If $\psi$ is an OSF for $\mathcal F$
          output by  \textsc{OSF-Builder}, then it is an sOSF.
	\end{itemize}
\end{lemma}
\pf
(i) We need to show that (P3) holds, that is, if $u \in V(G)$
is such that $\psi(u) \in V(T)$ for some
$T \in F$, then there exists some leaf $v \in \mathcal C_G(u)$ such that
$\phi(v) \in \mathcal C_T(\psi(u))$. But this clearly holds
since if $u \in V(G)$, by (S3),
we can take a child $v$ of $u$ in $G$ that is mapped by $\psi$
to a vertex below $\psi(u)$,
and then repeat this process of applying (S3) for $v$
and its children, until we reach a leaf of $G$.

(ii) In the top-down phase of \textsc{OSF-Builder}, at
a vertex $v$ we will always choose 
one of the children $u$ of $v$ to map to the same tree as $v$. Then apply (P2)
to see that $\psi(u)$ is a vertex below $\psi(v)$.
\epf

Before proving the main result of this section
(Theorem~\ref{th:minimum}), we give a useful characterization 
of the subsets of $A(G)$ of $G$ that can give rise to sOSF's.
For a set $I \subseteq A(G)$, we let $G-I$ denote the 
graph obtained by removing all arcs in $I$ from $A(G)$, 
which is a set of (not necessarily phylogenetic)  trees.
Suppose we are given a forest triple $\mathcal F = (G,F,\phi)$.
We call a set $I \subseteq A(G)$ an {\em introgression set
  for $\mathcal F$} if
the following hold: (i) if $u$ is a tail of some arc in $I$,
then $u$ is also a tail 
of some arc in $A(G)$ that is not in $I$, (ii) 
if $M$ is a tree in $G-I$, then $\phi$ maps every
leaf of $L(G)$ contained in $V(M)$
to the same tree  $T_M$ in $F$, and (iii) if
$M \neq M' \in  G-I$ and $u \in V(M)$, $v \in V(M')$ with $(u,v)$ or $(v,u)$  
in $I$, then $T_M \neq T_{M'}$. For example, the set $\{(w_1,a), (w_2,c)\}$
  is an introgression set for the forest triple depicted in
  Figure~\ref{fig:osf-not-sosf}. Note that analogous 
  sets for gene/species tree reconciliations (called transfer sets) are
 define in  \cite{tofigh2011simultaneous}.

It is straight-forward to 
check that for any sOSF $\psi$ for $\mathcal F$, the set of
arcs in $A(G)$ which map under $\psi$ to $C(\psi)$
is an introgression set. The following 
proposition which is analogous 
to \cite[Lemmas 4 and 5]{tofigh2011simultaneous}
for the DTL-scenario
shows that the converse holds. It  
implies that introgression sets give rise to a partition of 
the set of sOSF's for $\mathcal F$, where two 
sOSF's are in the same part if and only if they give rise
to the same introgression set for $\mathcal F$.

\begin{proposition}
	Suppose $\mathcal F = (G,F,\phi)$ is a forest triple 
	and $I \subseteq A(G)$. 	
	Then $I$ is an introgression set for $\mathcal F$
	if and only if there exists an sOSF $\psi$ for $\mathcal F$
	such that $C(\psi)$ coincides with the
          multiset $\{  (\psi(u),\psi(v)) \,:\, (u,v) \in I \}$.
\end{proposition}
\pf The only if statement follows from the remark preceding the proposition.
To see the converse, note that 
if $I$ is an introgression set for $\mathcal F$,
and we let $\psi_{I}:V(G) \to V(F)$ 
be the map given by setting, for all $u\in V(G)$,
$$
\psi_I(u) = lca_{T_M}(\{\phi(g)  \in L(T_M)
\,:\, g\in\mathcal C_G(u)\cap V(M)\},
$$
where $M$ is the tree in $G-I$ containing $u$. Then
$\psi_I$ is an sOSF.
\epf


As mentioned above, for a forest triple $\mathcal F$ there {exist
OSFs for $\mathcal F$ 
that are not sOSFs. Moreover, there exist sOSFs $\psi$ for
$\mathcal F$ with $|C(\psi)|=t(\mathcal F)$ which \textsc{OSF-Builder}
is not able to construct 
(\cite[Fig. 1, Supp. Mat.]{scholz2019osf}).
Even so, we now show that, for any forest triple,
\textsc{OSF-Builder} is guaranteed 
to produce an optimal OSF that is also an sOSF. 

\begin{theorem}\label{th:minimum}
	Suppose $\mathcal F=(G,F,\phi)$ is a forest triple. 
	Then 
	\begin{equation}
	\label{minimum}
	t(\mathcal F) = min\{|C(\psi)|  \,:\, \psi \mbox{ is a {\rm strict} OSF for } \mathcal F  \}.
	\end{equation}
	Moreover, \textsc{OSF-Builder} constructs an sOSF $\psi$ for $\mathcal F$
	with $|C(\psi)|=t(\mathcal F)$.
\end{theorem}
\pf
We begin by showing that Equation~(\ref{minimum}) holds. 
Given an OSF $\psi$ for $\mathcal F$, we put 
$$
U(\psi) = \{ u \in V^0(G) \,:\,  (\psi(u),\psi(v)) \in C^*(\psi)
\mbox{ for every child $v$ of $u$}  \}.
$$
Note that $\psi$ is an sOSF if and only if $U(\psi)= \emptyset$.
We claim that if $\psi$ is an OSF for $\mathcal F$ that is
not an sOSF, then there exists an OSF $\psi'$ for $\mathcal F$
such that $|C(\psi)| \ge  |C(\psi')|$ and  $|U(\psi)| > |U(\psi')|$.
Equation~(\ref{minimum}) then follows
since for any OSF $\psi$  for $\mathcal F$ we can keep applying
this claim until we obtain an sOSF $\psi''$ for $\mathcal F$
with $|C(\psi)| \ge |C(\psi'')|$ and $|U(\psi'')|=0$. 

To prove that claim, suppose that
$\psi$ is an OSF for $\mathcal F$ that is not an sOSF.
Choose some $u \in U(\psi)$ such that no vertex below $u$ but
distinct from $u$ is contained in $U(\psi)$. 
Then there exists a subset $\{u_1\dots,u_k\}$, $k\ge 1$,
of children of $u$ whose elements are all mapped by $\psi$ to 
some tree $T \in F$ which is different from 
the tree in $F$ containing $\psi(u)$. Define the map
$\psi_u:V(G)\to V(F)$ by setting $\psi_u(w)=\psi(w)$
for all $w \in V(G)-\{u\}$ and
$\psi_u(u)=lca_T(\psi(u_1),\dots,\psi(u_k))$. 
Then $\psi_u$ is an OSF for $\mathcal F$.
Note that this might have rendered the parent of $u$ 
an element in $U(\psi_u)$ so that $|U(\psi_u)|=|U(\psi)|$ holds. We
put $\psi=\psi_u$ and apply this construction of $\psi_u$  to a vertex 
$u\in U(\psi)$ such that no vertex below but distinct from
$u$ is contained in $U(\psi)$ and so on. Since $G$ is finite this implies that there must exist some OSF $\psi$ for $\mathcal F$ and
some $u\in U(\psi)$ such that $\psi'=\psi_u$ is an OSF for $\mathcal F$
and $|C(\psi)| \ge  |C(\psi')|$ and  $|U(\psi)| > |U(\psi')|$ holds
since, eventually, we will reach the root of $G$ which does 
	not have a parent. This  concludes the proof of  the claim.

Now, to see that the
second statement in the theorem holds, note that 
by Equality~(\ref{minimum}) it follows that if $\psi$
is any sOSF for $\mathcal F$
with $|C(\psi)|$ minimum, then $|C(\psi)| = t(\mathcal F)$. Hence, 
to complete the proof it suffices to show that
if $\psi'$ is any  sOSF for $\mathcal F$ output 
by \textsc{OSF-Builder}, then $|C(\psi)|=|C(\psi')|$. 

To see this, assume that $\psi$ is an sOSF for $\mathcal F$
such that  $|C(\psi)|$ minimum and that $\psi'$ is a
sOSF for $\mathcal F$ returned by \textsc{OSF-Builder}.
Then
clearly $|C(\psi')| \ge |C(\psi)|$.  To 
show that the reverse inequality
holds, first note that any introgression set $I$ in $A(G)$ gives rise in
a natural way to an extension $\bar{f}:V(G)\to F$ of the character 
$f_{\mathcal F}:L(G)\to F$. Indeed, we just extend $f_{\mathcal F}$
to $V(G)$ by taking the value of $v\in V(G)$
to be equal to that of $f_{\mathcal F}$ on the leaves of the tree in 
$G-I$ which contains $v$ in its vertex set. Hence, 
$|C(\psi)|\geq l_{f_{\mathcal F}}(G)$. But 
$l_{f_{\mathcal F}}(G) =|C(\psi')|$ since $\psi'$ is constructed by 
\textsc{OSF-Builder}. Thus, $|C(\psi)|\geq |C(\psi')|$ holds too
which implies $|C(\psi)|=|C(\psi')|$.
\epf

\section{Valid networks}
\label{sec:valid}

In applications, it is important to visualize OSFs using networks
as in Figure~\ref{fig:introduce}(c) in the introduction (see \cite{scholz2019osf}). 
In this section, we consider properties of networks that arise in 
this way from OSFs.

We begin by defining the network associated to an OSF.
Given a forest triple $\mathcal F=(G,F,\phi)$ and an OSF $\psi$
for $\mathcal F$, 
we define the network $N(\psi)$ to be the graph whose vertex set is $V(F)$
and whose arc set is $A(F)\cup C^*(\psi)$. An arc in $C^*(\psi)$
is called a {\em contact arc of  $N(\psi)$}. We illustrate these
concepts in Figure~\ref{fig:osf-not-sosf}. Note that 
the network of an OSF  may contain no roots and/or no leaves 
and it may also contain directed cycles (see e.g. Figure~\ref{fig:partrep}}).
In practice, this can make it tricky to interpret networks
but this can be alleviated by resolving them 
in a special way (see Appendix for more details).

We now present some 
basic properties of the network $N(\psi)$.


\begin{lemma}
  If $\psi$ is an OSF for a forest triple
  $(G,F,\phi)$, then $N(\psi)$ is a network on $L(F)$.	
Moreover, $N(\psi)$ is connected if and only if
$\phi(L(G)) \cap L(T) \neq \emptyset$ for all $T \in F$.
\end{lemma}
\pf
To show that $N(\psi)$ is a network, first note that Property~(N1) clearly 
holds, and that Property~(N2) holds because the trees in $F$
are phylogenetic trees whose leaf-set union is $L(F)$. 
Property~(N3) holds since every vertex in $L(F)$ is the head of some
arc in $F$ (as a
phylogenetic tree has at least 2 leaves).  Property~(N4) holds 
since clearly $\{ v \in V(N(\psi)) \,:\, outdeg(v)=0 \} \subseteq  L(F)$.

We now show that the second statement holds. Suppose first 
that $\phi(L(G)) \cap L(T) \neq \emptyset$ for all $T \in F$.
By Property~(P2), every arc $(u,v)$ in $G$ is either
mapped under $\psi$ to an element in $C^*(\psi)$
or $\psi(u)$ and $\psi(v)$ are both contained in the same tree $T$ of $F$
with $\psi(v)$ below $\psi(u)$ in $T$. Hence, if for
any tree in $F$ we pick some leaf $v\in \phi(L(G))$
(which is possible as $\phi(L(G)) \cap L(T) \neq \emptyset$ for
  all trees $T\in F$), it
follows that under $\psi$ the path in $G$ from 
$\rho_G$ to any $u \in L(G)$ with $\phi(u)= v$ yields
a path in $N(\psi)$ from $\psi(\rho_G)$ to $v$.
In particular, for each tree $T$ in $F$, there is an
undirected path in the underlying 
graph of $N(\psi)$  from $\rho_G$ to some 
vertex in $T$ in $F$. Hence, $N(\psi)$ is connected.

Conversely, suppose that $N(\psi)$ is connected.
Let $T\in F$. If $\psi(\rho_G) \in V(T)$ then,
by Property~(P3), there is some $v \in \mathcal C_G(\rho_G)=L(G)$ such that
$\phi(v)=\psi(v) \in \mathcal C_T(\psi(\rho_G))$. So 
$\phi(L(G)) \cap L(T) \neq \emptyset$.
So assume that $\psi(\rho_G) \not\in V(T)$. Let $T'$
denote the tree in $F$ that contains $\psi(\rho_G)$ in its vertex set.
As $N(\psi)$ is connected and $|F|\geq 2$, the construction of $N(\psi)$
implies that there must exist a path from $\psi(\rho_G)$
  to a vertex $v$ in $T$. Without loss of generality we may assume that
  $v$ is the head of some arc in $C^*(\psi)$. Then, by the definition
  of a contact arc, there must exist some
  arc $(w,u)$ in $G$ such that $\psi(u)=v$. Since $\psi(u)\in V(T)$,
  Property~(P3) implies that there must exist some leaf $x\in \mathcal C_G(u)$
  such that $\phi(x)\in \mathcal C_T(\psi(u))$. Hence, 
  $\phi(L(G)) \cap L(T) \neq \emptyset$ must hold in this case too.
\epf

Now, we say that a network $N$ is {\em a representation} of an OSF
$\psi$ on $(G,F,\phi)$ if it is isomorphic to $N(\psi)$ or -- equivalently --  
there is a set $A \subseteq A(N)$ such that $N-A=(V(N),A(N)-A)$ is 
isomorphic to $F$ and $A = C^*(\psi)$ (under the isomorphism 
between $N-A$ and $F$).  
In addition, we call a network $N$ {\em valid} if there exists a forest
triple $\mathcal F$ such that $N$ is a  
representation of some OSF on $\mathcal F$.

Note that a network can be a representation for more than one OSF
(see e.g. Figure~\ref{fig:partrep}). Therefore, it is interesting
to know which types of OSF 
can be the representation of a valid network.
For example, note that a valid network is always a 
representation of some sOSF. 
Indeed, suppose $N$ is of the form $N(\psi)$ for 
some OSF $\psi$ for some forest triple $(G,F,\phi)$. Then we can
create a new forest triple $(G',F,\phi')$
by inserting a new leaf vertex $g_u$ in $G$ pendant to each $u \in V^0(G)$ 
to create $G'$ and extending $\phi$
to a leaf-map on $L(G')$ by mapping 
each new vertex $g_u$ to a leaf in the tree to which $u$ is mapped under $\psi$ that is below $\psi(u)$. 
Then it is straight-forward to check that the map $\psi'$ obtained
by extending $\psi$ in the natural way to
a map $V(G')\to V(F)$ is an sOSF for $(G',F,\phi')$
such that $N(\psi)$ is isomorphic to $N(\psi')$. 

Interestingly, as we shall now show, an even stronger result holds:

\begin{theorem}\label{thm:valid}
  Suppose that $\psi$ is an OSF for some forest triple
  $\mathcal F=(G,F,\phi)$. Then
  there is an sOSF $\psi'$ for some forest triple $\mathcal F'=(G',F,\phi')$
  that is output by \textsc{OSF-Builder} such that
  $N(\psi)$ is isomorphic to $N(\psi')$. In particular, a network $N$ is valid
  if and only if there is some OSF $\psi$
  output by \textsc{OSF-Builder} such that $N$ is
  a representation of $\psi$. 
\end{theorem}
\pf
Consider the tree $G$. 
Construct a new phylogenetic tree $G'$ by inserting,
for each $u \in V^0(G)$, 
$outdeg(u)+1$ new arcs into $G$ 
of the form $(u,v)$ (so that in particular $v$ is a leaf in $G'$). 
Extend the leaf-map $\phi$ on $L(G)$ to a leaf-map
$\phi'$ on $L(G')$ by,
for each $u$ in $V^0(G)$,  mapping the new children 
of $u$ arbitrarily onto a set $S$ of leaves in $L(T)$, where 
$T \in F$ is the tree with $\psi(u) \in V(T)$, so that $lca_T(S)=\psi(u)$.
Let $\psi'$ be the extension of $\psi$ to $V(G')$ given by putting,
  for all $v\in V(G')$,  $\psi'(v)=\psi(v)$ if $v\in V(G)$
and $\psi'(v)=\phi'(v)$ otherwise. 

It is straight-forward to check that 
$\psi'$ is an sOSF for $(G',F,\phi')$, and that 
$N(\psi')$ is isomorphic to $N(\psi)$. 
Moreover, for all $u \in V^0(G')$, the definition
of $G'$ implies that  the set $\sigma(u)$  computed 
by \textsc{OSF-Builder} must have size 1 and that  it 
consists of the tree in $F$ in which $\psi(u)$ is contained.
It follows that the necessarily
unique map obtained by \textsc{OSF-Builder}
by applying its top-down phase  is equal to $\psi'$. 
\epf

\section{A characterization of valid networks}  
\label{sec:representing}

In this section we give a characterization for valid networks.
To do this we first show that a valid network can be
represented by a special type of OSF.

We start with a definition. Suppose that $(G,F,\phi)$ is a forest triple.
Given a path $\gamma = w_1, w_2,\dots, w_k$, $k\geq 2$, in $G$ we
obtain a walk $\gamma'$ in $N(\psi)$ (possibly of length 0) 
that starts at $\psi(w_1)$ and ends at $\psi(w_k)$. 
Central to our definition is the observation that if
$(w_i,w_{i+1})$, $1\leq i\leq k-1$, is a contact arc in $G$, then 
$(\psi(w_i),\psi(w_{i+1}))$ is an arc in $N(\psi)$ but not in $F$.
Otherwise $\psi(w_i)$ and $\psi(w_{i+1})$ are both vertices in
some tree $T$ of $F$, 
and, therefore, there exists a path in $T$ (possibly of length 0) 
starting at $\psi(w_i)$ and ending at $\psi(w_{i+1})$. 
To obtain $\gamma'$ we then take the walk obtained
by inserting these paths into the sequence
$\psi(w_1),\psi(w_2),\dots,\psi(w_k)$, where any 
consecutive repeats are suppressed. 

\begin{proposition}\label{trail}
  Suppose that $\psi$ is an OSF for some forest triple
  $(G,F,\phi)$. Then there is 
  an OSF $\psi'$ for some forest triple $(G',F,\phi')$ such that
  $\gamma'$ is a trail in $N(\psi')$
  for every path $\gamma$ in $G$, and $N(\psi)$
  is isomorphic to $N(\psi')$.
\end{proposition}
\pf
To help construct the forest triple $(G',F,\phi')$,
  we first employ Property~(P3) to choose
for every vertex $v\in V(G)$ such that $\psi(v)$
  is a non-leaf vertex of some tree $T \in F$
some leaf $l_v\in L(G)$ below $v$ such that
$\phi(l_v)\in\mathcal C_T(\psi(v))$.
 Subsequently, we insert
an arc $(v,u)$ into $G$ for each leaf in $T$ that is below $\psi(v)$ in
$T$ and denote the resulting tree by $G'$. Note that
if $\psi(v)\not\in L(F)$ then $outdeg_{G'}(v)\geq 4$.
By extending $X=L(G)$ canonically we may assume that
$G'$ is a phylogenetic tree.

Next, we define a new leaf-map $\phi': L(G')\to L(F)$
  by putting, for all $l\in L(G')$, 
  $\phi'(l)=\phi(l)$ if $l\in L(G)$ and $\phi'(l)=\phi(l_v)$
  otherwise, where $v\in V^0(G)$ is the vertex
  to which the arc $(v,l)$ was attached in
  the construction of the extended tree $G'$.  By construction
it follows that $(G',F,\phi')$ is a forest triple and that
$\psi$ can be canonically extended to an OSF $\psi'$ for $(G',F,\phi')$.
Note that $N(\psi)$ and $N(\psi')$ are clearly isomorphic
by construction of $G'$. To reduce notation, 
we shall now denote $G'$,  $\phi'$ and $\psi'$
by $G$, $\phi$ and $\psi$, respectively.

Suppose now that there is some  path $\gamma$ in $G$
for which the associated walk
$\gamma'$ in $N(\psi)$ is not a trail. Then there must 
exist some minimum length
subpath $\delta=v_1,\dots,v_{k+1}$, $k\geq 3$, of $\gamma$ 
that gives rise to a walk $\delta'$ in $N(\psi)$. Without loss of
  generality we may assume that $\delta$ is such that 
  $\psi(v_1)\not=\psi(v_2)$ and  $\psi(v_k)\not=\psi(v_{k+1})$,
  that the first  and last arc of $\delta'$ coincide, and that
  no arc on the induced walk  $v_2,\dots,v_k$ of $\delta'$ occurs
  more than once. Note that $\delta'$ must contain at least
  	four non-leaf vertices. Also note that
  in particular $(\psi(v_1),\psi(v_2))=(\psi(v_k),\psi(v_{k+1})) \in C^*(\psi)$
  as otherwise Property~(P2) would not hold for $v_2$ and $v_k$.

\begin{figure}[h]
\begin{center}
\includegraphics[scale=0.7]{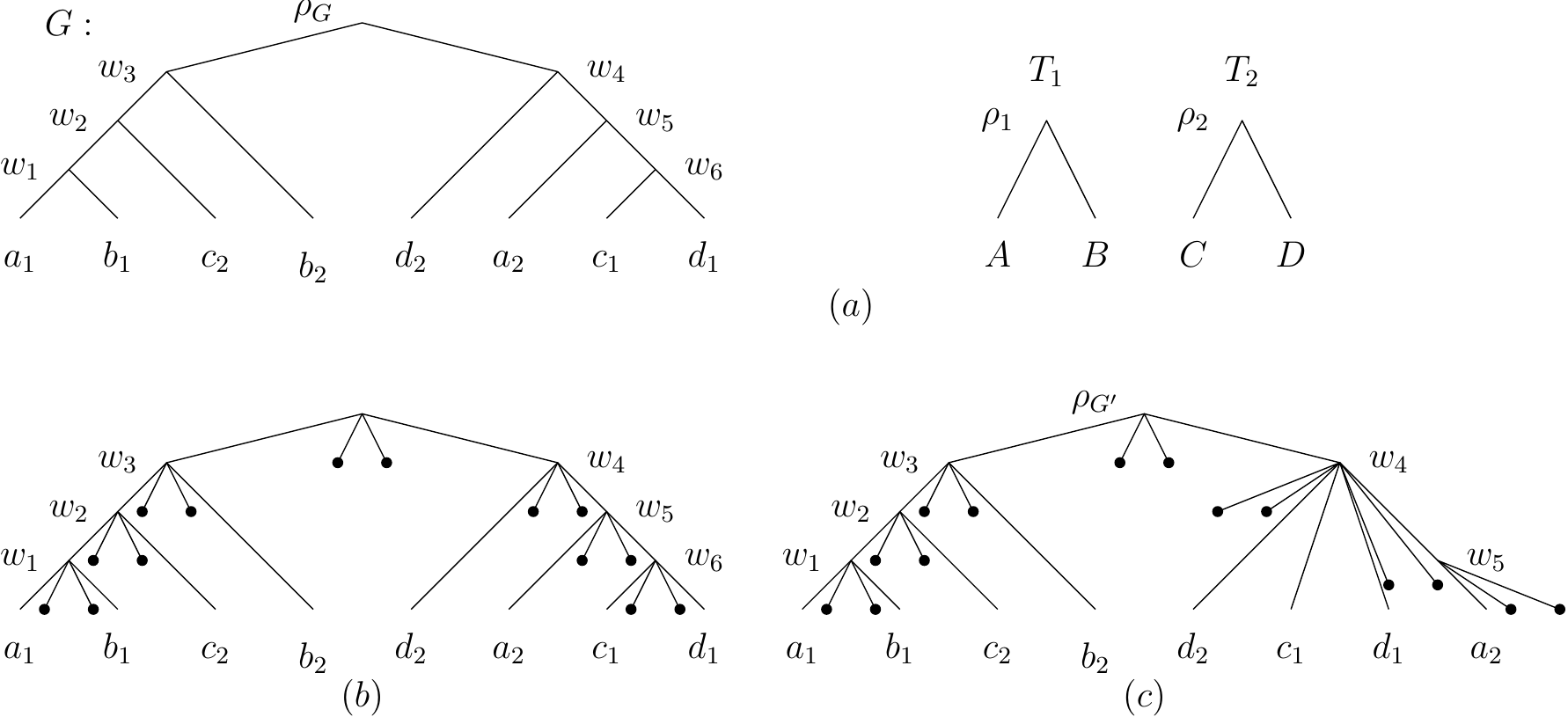}
\caption{
  (a) The map $\psi:V(G)\to V(F)$
 defined as $\psi(\rho_G)=\psi(w_i)=\rho_1$ for all $i=1,2,3$, and 
 	$\psi(w_i) =\rho_2$, for all $i=2,4,6$ is an OSF for the forest triple given by the depicted tree $G$, forest $F=\{T_1,T_2\}$, and
 	$\phi$ the usual leaf-map. (b) The phylogenetic tree
  obtained from $G$ in the first part of the construction in
  Proposition~\ref{trail}. (c) The phylogenetic tree $G'$ obtained
  from the tree in (b) in the second part that construction.}
\label{fig:construction}
\end{center}
\end{figure}
Now, remove the arc $a=(v_k,v_{k+1})$ from $G$, and
  if $v_{k+1}\in L(G)$ then also replace $a$ by $(v_2,v_{k+1})$
  and if $v_{k+1}\not\in L(G)$ then replace
each arc of the form $(v_{k+1},u)$ by the arc $(v_2,u)$. 
It follows that the resulting tree $G'$ is a phylogenetic tree.
We illustrate the construction of $G'$ in Figure~\ref{fig:construction}, where the forest $F$ is obtained by deleting the dashed arcs in the network displayed in Figure~\ref{fig:construction}(b) and the map
	$\psi:V(G)\to V(F)$ 
	given by $\psi(\rho_G)=\psi(w_i)=\rho_1$ for all $i=1,2,3$, and 
	$\psi(w_i) =\rho_2$, for all $i=2,4,6$, is an OSF for the
	forest triple $(G,F,\phi)$ where $G$ is the tree depicted
	in Figure~\ref{fig:construction}(a) and $\phi$ is the usual leaf-map.
	Note that the walk $\gamma'$ in $N(\psi)$
	associated to the path $\gamma:=\rho_G,w_4,w_5,w_6$ in $G$
	is not a trail since the arc $(\rho_1,\rho_2)$ is crossed twice.

Continuing with the proof, we next adjust $\phi$ to become a leaf-map $\phi': L(G')\to L(F)$ from 
$G'$ into $F$ by putting, for all $x\in L(G')$, $\phi'(x)=\phi(x)$
if $x\in L(G)$ and $\phi'(x)=\phi(l_v)$ otherwise,
where $v\in V(G)$ is the 
vertex to which the arc $(v,x)$ was attached in the construction of $G$,
we obtain a  forest triple $(G',F,\phi')$.
Since the map $\psi':V(G')\to V(F)$ given by putting,
  for all $v\in V(G')$, $\psi'(v)=\psi(v)$ if $v\in V(G)$
  and $\psi'(v)=\phi'(v)$ otherwise, clearly 
  satisfies Properties~(P1) --
  (P3) by construction and because $\psi$ is an OSF for $(G,F,\phi)$,
  it follows that $\psi'$ is an OSF for $(G',F,\phi')$. Note that, also
  by construction we must have that $N(\psi')$ is isomorphic to $N(\psi)$.


To conclude the proof, note that 
$|A(G')|\leq |A(G)|$, and, by construction, $\psi'$ maps
one less arc in $G'$ to the contact arc
$(\psi(v_1),\psi(v_2))$ of $N(\psi')$. 
Hence, as $G$ is finite, we can keep repeating this process until we eventually 
end up with a forest triple $\mathcal F$ and some OSF for $\mathcal F$
satisfying the conditions given in the statement of the  proposition.
\epf

\begin{remark}
  Note that in general, we cannot replace ``$\gamma'$ is a trail" by
  ``$\gamma'$ is a path" in
the statement of Proposition~\ref{trail}. This can be shown by considering 
the example depicted in Figure~\ref{fig:construction},
where it can be seen
that every possible OSF 
results in some path $\gamma$ for which $\gamma'$ is not a path.
\end{remark}

We now present the aforementioned characterization of valid networks.
Our proof uses the idea of unfolding a network, which was 
introduced in \cite{HM06}.

\begin{theorem}\label{thm:structure}
Suppose that $N$ is a network on $X$. Then $N$ is valid if and only if 
there is some vertex $\rho \in V(N)$  and some $A \subseteq A(N)$ such that 
\begin{itemize}
\item[(V1)] $N-A=(V(N),A(N)-A)$ is a forest $F$ with $|F|\ge 2$
  and leaf-set $X$
	in which every arc in $A$ has its ends in different trees of $F$, and
      \item[(V2)] every arc in $A$ is contained in some trail
        $\rho=v_1,v_2,\dots,v_m$, $m \ge 2$,
	such that if $v_i,v_j$, $i \le j$ are vertices in the same
        tree $T$ in $F$	then $v_i \preceq_T v_j$.
\end{itemize}
\end{theorem}

\pf  Suppose first that $N$ is valid, that is,
$N=N(\psi)$ for some OSF $\psi$ of a forest triple $\mathcal F = (G,F,\phi)$. 
Put $A=C^*(\psi)$. As $\psi$ is an OSF for $\mathcal F$, it follows
in view of Property~(N3) that
$N-A$ is a forest $F$ such that
	every arc in $A$ has its ends in different trees of $F$. 
In view of Properties~(N2) and (N4), we also have
  $L(F)=X$.  Hence, (V1) holds.

We now show that (V2) holds
relative the vertex $\rho =\psi(\rho_G)$ and the set $A=C^*(\psi)$. 
Note first that by Proposition~\ref{trail}, 
we may assume 
for every path $\gamma$ in $G$ that the associated walk
$\gamma'$  in $N(\psi)$ is a trail.
Now, suppose we are given 
some arc $(u',v')$ in $A$. Then by definition of $C^*(\psi)$, there
must exist some arc $(u,v)$ in $G$ such that
$\psi(u)=u'$ and $\psi(v)=v'$. Let $T$ denote the tree in $F$ that
  contains $v'$ in its vertex set. Then Property~(P3) implies that
  there must exist some leaf $l\in \mathcal C_G(v)$ such that
  $\phi(l)\in\mathcal C_T(\psi(v))$.  Hence, there must exist a
  path $\gamma$ in $G$  from $\rho$ to $l$ that
  contains $(u,v)$. By (P2) and the first observation
  in this proof, it follows that the trail
  $\gamma'$ in $N(\psi)$  associated to $\gamma$ has the 
required properties. So (V2) holds, as required.

Conversely, suppose $N$ is a network on $X$, that 
contains an arc set $A$ and a vertex $\rho$ 
such that (V1) and (V2) are satisfied. Let $F$ denote
the forest $N-A$. Then $|F|\geq 2$ and $L(F)=X$.

We start with constructing a forest triple $(G,F,\phi)$.
Let $G=G(N)$ be the graph obtained by
``unfolding $N$ at $\rho$", which 
we define as follows (cf. \cite[p.617]{HM06}):

\begin{itemize}
	\item the vertices of $G$ are trails 
	$\rho=v_1,v_2,\dots,v_m$, $m \ge 1$,
	  in $N$ such that if $v_i,v_j$, $i \le j$, are both
          vertices in some tree $T$ of $F$
	then $v_i \preceq_T v_j$;
      \item for all vertices $\gamma$ and $\gamma''$ in $G$
        there is an arc from $\gamma$ to $\gamma''$ in $G$ 
	if and only if $\gamma'' =  \gamma a$
	holds for some arc $a$ in $A(N)$;
	\item the vertices in $G$ that start at $\rho$ and end at a
	vertex $v\in X$ are labelled by distinct elements in 
	the set $X_v = \{l_i \in L(N)\,:\, 1 \le i \le n_v\}$, 
	where $n_v$ is the 
	number of vertices in $G$ which end at $v$.
\end{itemize}

By construction it immediately follows that $G$ is a tree with root 
$\rho$ (considered as a path with length 0). 

To see that $G$ is a phylogenetic tree suppose that $\gamma$
is an element in $V(G)$ whose end vertex $v$ 
is not a leaf in $N$. Then, by (V1), $v$ must be an interior
vertex of some tree $T$ in $F$.  Hence there must 
be at least 2 possibilities to extend $\gamma$ by adding an arc
in $T$ whose tail is $v$. Hence, $\gamma$ is an interior vertex of
$G$ with outdegree at least 2. 
Moreover, the leaf-set of $G$ is equal to 
the set obtained by taking the union $Y$ of the sets $X_v$,
where the union is taken over
all $v\in X$
Since $|Y|\geq 2$, it follows that $G$ is a phylogenetic tree on $Y$. 

To obtain the leaf-map $\phi:L(G)\to L(F)$ note that
there is a natural map from $V(G)$ into $V(N)$, which 
takes each element in $V(G)$ to its end vertex in $N$. This
gives rise to a map $\psi$ from $V(G)$ to $V(F)$,
which maps the vertex set of $G$ to the vertex set of $N$.
We let $\phi = \psi |_{L(G)}$ be the map from $L(G)$ to $L(F)$
induced by $\psi$. This completes the  
construction of the forest triple $(G,F,\phi)$.

We now claim that $\psi$ is an OSF for $(G,F,\phi)$ such that $A = C^*(\psi)$, 
which will complete the proof the theorem.
Property~(P1) holds by definition of $\phi$. 
The fact that Property~(P2) holds follows immediately from
(V2) and the definition of $G$ and $\psi$. 
Property~(P3) holds since if $\gamma$ is an element 
in $V(G)$ with end vertex $v$ in $N$, then $v$ is contained
in some tree $T$ of $F$. In case $v$ is not a leaf of $T$ 
there must exist a leaf $u\in V(T)$ below $v$.
In that case, let $\gamma_1$ in $V(G)$ denote the path
obtained by extending $\gamma$ by some path in $T$ from $v$ to $u$ and put $\gamma_1=\gamma$, otherwise. In particular, 
it follows that $\gamma_1\in {\mathcal C}_G(\gamma)$ 
and $\phi(\gamma_1)= \psi(\gamma_1) \in {\mathcal C}_T(v)$.

It remains to show that $A = C^*(\psi)$.
Clearly $C^*(\psi) \subseteq A$ by (V1).
Conversely, suppose $a \in A$. By (V2), $a$ 
is contained in some vertex $\gamma$ in $V(G)$, which 
we may assume without loss of generality to have $a$ as the last arc.
But then for the path $\gamma_1$ in $V(G)$ with $\gamma = \gamma_1 a$
we have $(\gamma_1,\gamma) \in A(G)$ and
$a = (\psi(\gamma_1),\psi(\gamma))$. Thus $a \in C^*(\psi)$.
\epf

\section{Stability of optimal scores under tree and forest alterations}
\label{OSF:stab}

In this section we are motivated by the following question. If 
we alter $G$ or $F$ in
a forest triple $\mathcal F = (G,F,\phi)$, then how much does
this change the minimum number of contact arcs of an
OSF for $\mathcal F$ (i.e. by how much can $t({\mathcal F})$ change)?
This is important since it can help indicate, for example, how
changes to the input of \textsc{OSF-Builder} can effect its output.   


From now on all forest triples are assumed to be binary.
We begin by recalling some facts concerning 
subtree prune and regraft (SPR) operations for 
rooted phylogenetic trees as defined in \cite{BS04}. 
An {\em SPR operation} on a tree $T$ is defined by
cutting any arc $(u,v)$ in $A(T)$ (so pruning off the 
subtree of $T$ rooted at $v$), and then regrafting this pruned subtree
into a subdivided arc of the pruned tree. Note that some
care has to be taken when pruning off a subtree adjacent to
the root; the reader can find full details in \cite{BS04}.
The {\em rooted SPR distance} between two phylogenetic
trees $T$ and $T'$ on the same leaf set, denoted $d_{rSPR}(T,T')$, is
the minimum number of SPR operations requited to transform
one tree into the other.

We now ask the above question more precisely: If
$G$ and $G'$ are two phylogenetic trees that
differ by $k$ SPR operations and
$\mathcal F = (G,F,\phi)$ and $\mathcal F' = (G',F,\phi)$ are
forest triples, then how different
can $t(\mathcal F)$ and $t(\mathcal F')$ be in terms of $k$?
We next give an upper bound for this difference:

\begin{theorem}\label{spr-tree}
	If $G$ and $G'$ are two binary phylogenetic trees on the same leaf-set, and
	$\mathcal F = (G,F,\phi)$ and $\mathcal F' = (G',F,\phi)$ are
	forest triples, then 
	$|t(\mathcal F)-t(\mathcal F')| \le d_{rSPR}(G,G')$.
\end{theorem}

This theorem follows immediately from \cite[Corollary 3.12]{FK16}. This
states that if $G$ and $G'$ are two binary phylogenetic trees on $X$,
and $f$ is a character on $X$, then $|l_f(G)-l_f(G')|\le d_{rSPR}(G,G')$.
The theorem thus follows since if $\mathcal F$ is a forest 
triple, then $t(\mathcal F) = l_{f_\mathcal F}(G)$. 

The results in \cite{FK16} also provide the following additional 
bounds on $|t(\mathcal F)-t(\mathcal F')|$. In case 
there are $r$ trees in $F$ whose leaf-sets intersect the image $\phi(L(G))$, 
and $r \le n=|L(G)|$,
then by \cite[Lemma 3.14]{FK16} 
$$
|t(\mathcal F)-t(\mathcal F')| \le \lfloor (r-1)(\frac{n}{r} -1)\rfloor,
$$
and by \cite[Theorem 3.15]{FK16}  
$$
|t(\mathcal F)-t(\mathcal F')| \le n- 2\sqrt{n} +1,
$$
a bound which is tight for $n=9$ (\cite[Figure 4]{FK16}). 

We now present a result which describes
how applying an SPR operation to $F$ in a forest 
triple $\mathcal F = (G,F,\phi)$ can affect $t({\mathcal F})$.  
Note that if we alter any tree in $F$ by a SPR operation
then this has no effect on 
$t({\mathcal F})$, as the associated character $f_{\mathcal F}$ remains unchanged.

\begin{theorem}\label{spr-forest}
	Let $\mathcal F = (G,F,\phi)$ be a forest triple.
	Suppose a forest $F'$ is obtained from $F$ by 
	pruning off some subtree $T_0$ in a tree $T \in F$ and grafting
	$T_0$ into a tree in $F - \{T\}$. Let
	$\mathcal F' = (G,F',\phi)$ denote the resulting forest triple. Then 
	$$
	|t(\mathcal F)-t(\mathcal F')| \le  |\{ v \in L(G) \,:\, \phi(g) \in L(T_0)\}|.
	$$
\end{theorem}

The theorem follows immediately from the following lemma, which 
generalizes \cite[Observation 4.2]{FK16}.

\begin{lemma}
	If $f$ is a character on a set $X$, $T$ is a phylogenetic 
	tree on $X$, and $f'$ is a character on $X$ obtained by
	changing the value of $f$ on exactly $k \ge 1$ elements in $X$, 
	then $|l_{f'}(T)-l_f(T)| \le k$.
\end{lemma}
\pf
We use induction on $k$. The base case, $k=1$, is
\cite[Observation 4.2]{FK16}.

Now, suppose the inequality in the lemma holds for 
all $k \le L-1$, some $L \ge 2$. Let $f$ be a character on
$X$, let $T$ be a phylogenetic 
tree on $X$, and let $f'$ be a character on $X$ obtained by
changing the value of $f$ on exactly $L$ elements in $X$.
Let $f''$ be a character on $X$ obtained by
changing the value of $f'$ on exactly one element $x \in X$
to the value $f(x)$. Then $|l_{f''}(T)-l_{f'}(T)| \le 1$,
and by induction $|l_{f'}(T)-l_{f}(T)| \le L-1$. The lemma follows.
\epf

\begin{remark}
  The bounds obtained in Theorems~\ref{spr-tree} and \ref{spr-forest}
  are  in agreement with the 
	results obtained in a 
	simulation study in \cite{scholz2019osf}, where the affect of 
	altering $G$ and $F$ by SPR operations on $t((G,F,\phi))$
        was studied empirically using simulations.
\end{remark}

\section{Conclusion}
\label{sec:conclusion}

There are several new directions that could be of potential interest. 
One possibility following on from the results 
presented in the last section could be to try and understand the effect 
that changing the root location of the gene tree can have on 
the optimal score for an OSF for a forest triple. This is important as it can
be difficult to accurately root the tree in practice, and 
its location is known to have an impact on duplication-transfer-loss models for
lateral gene transfer (see e.g.\cite{kundu2018impact}).

In another direction, it could be worth defining and
studying spaces of OSF's for a given forest triple -- such spaces
have been intensively studied for gene/species tree reconciliation models and 
can, for example, give important insights on how optimal
OSF's are distributed over the collection of 
all possible OSFs  (see e.g. \cite{doyon2009space}).
Defining such spaces
could also lead to new metrics that could be used to 
compare OSF's such as the ones defined in e.g. 
\cite{huber2019exploring}.

 Finally, note that a network for an OSF can 
be thought of a special example of a forest with extra
edges added. There has been much work recently
on understanding the structure of networks that arise
from adding some edges into a tree (so-called tree-based networks) 
\cite{francis2015phylogenetic,zhang2016tree}.
It would thus be of interest to see what properties ``forest-based" networks
have in common (or not) with tree-based networks.\\

\section*{Appendix}\label{sec:appendix}

\begin{figure}[ht!]
	\center
	\scalebox{0.7}{\includegraphics{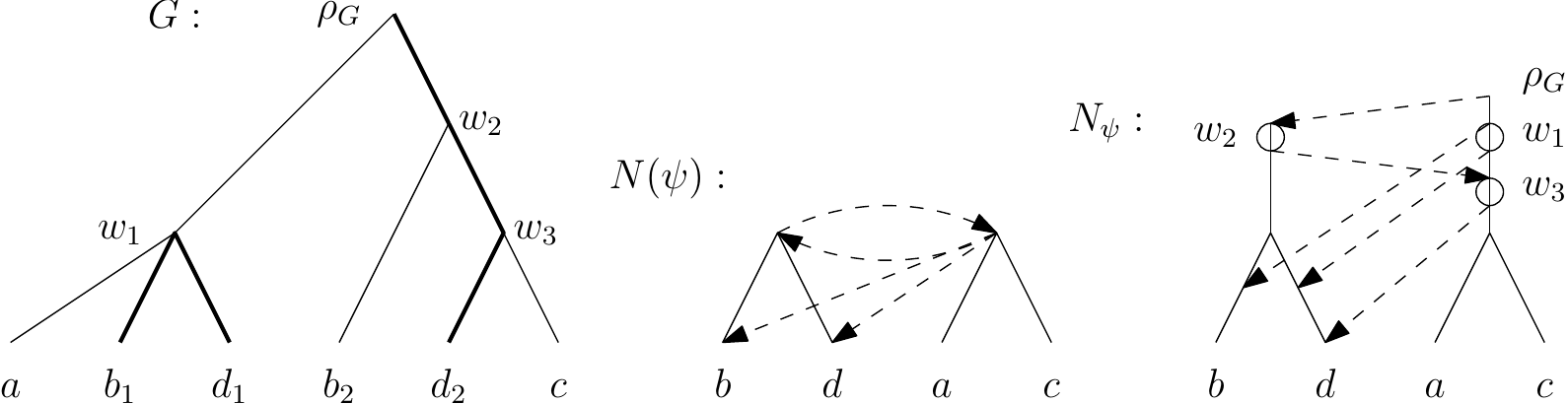}}
	\caption{The networks $N(\psi)$ and $N_{\psi}$
            for the OSF $\psi$ given in the text
            where $G$ is the depicted tree,  $F$ is the forest
			obtained by deleting the dashed arcs in $N(\psi)$ and
			$\phi$ is the usual leaf-map. 
			The four dashed	arcs make up $C^*(\psi)$ and the
                        five edges in bold in $G$ are the induced
                        introgression set. Each circle in $N_{\psi}$ indicates 
                        a set of subdivision vertices all of which are labelled by the vertex next to it.
                        }	
	\label{fig:incidental}
\end{figure}

As mentioned in Section~\ref{sec:valid}, the network $N(\psi)$ for an OSF  $\psi$ 
may have no roots and/or no leaves and it may also contain directed cycles
which can make it difficult to interpret in practice.
We now briefly 
present a way to produce a binary resolution $N_{\psi}$
of $N(\psi)$ which can help to circumvent this issue (cf. \cite{scholz2019osf}). 
We also show that this resolution has the attractive
property that it only contains {\em incidental cycles} in $N(\psi)$,  
that is, cycles that are {\em not} of the form $\gamma'$ for some 
path $\gamma$ in $G$ (see e.g. Figure~\ref{fig:incidental}, 
where for the OSF $\psi:V(G)\to V(F)$ given by mapping $\rho_G$, $w_1$ and $w_3$  to the parent of $a$ and $c$ in $N(\psi)$ and $w_2$ to the parent of
	$b$ and $d$ in that network, the walk $\psi(\rho_G), \psi(w_2), \psi(w_1)=\psi(\rho_G)$
is an incidental cycle in $N(\psi)$.
This is important, as it would not make evolutionary sense
for a gene in some species to introgress into one of its ancestors.

We begin with a simple observation that is implied by Property~(P2).

\begin{observation}\label{helper}
  Suppose that $\psi$ is an OSF for $\mathcal F=(G,F,\phi)$ and that $\gamma$
  is a path in $G$ such that $\gamma'$ is a cycle in  $N(\psi)$, then  the
  first and last arcs in $\gamma'$ are both contact arcs for $\psi$.
\end{observation}

We now explain how to generate a binary
  resolution $N_{\psi}$ of $N(\psi)$
  -- see Figure~\ref{fig:incidental} for an example.
	Suppose $\psi$ is an OSF for some forest triple
	$\mathcal F=(G,F,\phi)$. Without loss of generality we may assume
        that $\psi$ is an sOSF. Let $I\subseteq A(G)$
        denote the introgression set
        induced by $C(\psi)$ on $G$.
 We start with associating a directed graph $N'$ to $N(\psi)$.
        To do this, we first attach an
	incoming arc to every root of $N(\psi)$. We consider these arcs
	also as arcs in $F$. For every vertex $v\in V(N(\psi))$
        such that there exists some $w\in V(G)$ with $v=\psi(w)$
        we then add $|\{a\in I\,:\, w\in\{head(a),tail(a)\}\}|$
        subdivision vertices to the incoming arc of $v$ in $F$
        all of which we label $w$. If $v$ is in the image of more than
        one vertex of $G$ under $\psi$, say $w$ and $w'$,
        then we also ensure that the ancestral relationships
        between the subdivision vertices labelled $w$ and $w'$ are preserved.
        Next, we attach every  contact arc in
        $C^*(\psi)$ to a pair of correspondingly labelled subdivision
        vertices so that every arc in $G$ is an arc in the resulting
        graph. Finally, we remove all vertices of
        indegree zero and outdegree one plus their
	outgoing arcs. 

Note that the resulting graph $N'$
	has at least one root vertex, leaf set $L(F)$, and potentially still 
	vertices that are involved in two or
	more outgoing arcs in $F$. Also note that  $N'$
	might not be a representation of $\psi$ since the trees in 
	$N'-C^*(\psi)$ might contain vertices with indegree and outdegree one.
	However $N'$ can be easily transformed into a representation
	of an OSF $\psi'$ for some forest triple $(G,F',\phi')$
        with $C^*(\psi')=C^*(\psi)$ by attaching
	to each subdivision vertex of $N'$ a new leaf to obtain a
	new network $N''$ and defining $F'$, $\phi$ and $\psi'$
        in the canonical way.
%
        
Note that by Observation~\ref{helper}, any
	cycle in $N(\psi)$ that is {\em not} incidental does not give
	rise to a cycle in $N''$ (i.e. non-incidental cycles in $N(\psi)$ are
	broken in $N''$). Resolving
	potential vertices in $N(\psi')$ that have three or more
	outgoing arcs in $F$ results in a binary resolution 
	for $N(\psi')$. The directed graph obtained
	by removing all leaves from that resolution 
	that are not also contained in $N(\psi)$ (plus their incident
	arcs) is $N_{\psi}$.
%
%


In this
construction some choices might be made
and so 
there could be several binary resolutions. However, 
every binary resolution has the following property, which 
follows from the above construction.

\begin{proposition}
	Suppose that $\psi$ is an OSF. Then any cycle in $N_{\psi}$ 
	must come from an incidental cycle in $N(\psi)$.
\end{proposition}

This proposition  highlights  the way in which the OSF 
model differs from the model of lateral gene transfer 
presented in \cite{tofigh2011simultaneous}. In the latter 
model biologically unfeasible transfer scenarios can potentially 
arise \cite[cf. p.520]{tofigh2011simultaneous}  (essentially a gene 
that transfers into an ancestor of the species from which it arose).
However, by the last result this is not the case 
for the OSF model, since we only obtain incidental cycles in $N_{\psi}$. \\

\noindent{\bf Acknowledgement}
The authors would like to thank Martin Taylor and Andrei-Alin Popescu
for stimulating discussion on the topic of the paper. GS would like to
thank the School of Computing Sciences, University of East Anglia 
where some of the ideas underlying this paper were conceived.

\bibliographystyle{plain}
\bibliography{osf}

\begin{thebibliography}{10}

\bibitem{BS04}
M.~Bordewich and C.~Semple.
\newblock On the computational complexity fo the rooted subtree prune and
  regraft distance.
\newblock {\em Ann. Combin.}, 8:409--423, 2004.

\bibitem{doyon2009space}
J.-P. Doyon, C.~Chauve, and S.~Hamel.
\newblock Space of gene/species trees reconciliations and parsimonious models.
\newblock {\em Journal of Computational Biology}, 16(10):1399--1418, 2009.

\bibitem{FK16}
M.~Fischer and S.~Kelk.
\newblock On the maximum parsimony distance between phylogenetic trees.
\newblock {\em Ann. Combin.}, 20(1):87--113, 2016.

\bibitem{F71}
W.~Fitch.
\newblock Toward defining the course of evolution: minimum change for a
  specified tree topology.
\newblock {\em Syst. Zool.}, 20:406--416, 1971.

\bibitem{francis2015phylogenetic}
A.~R. Francis and M.~Steel.
\newblock Which phylogenetic networks are merely trees with additional arcs?
\newblock {\em Systematic Biology}, 64(5):768--777, 2015.

\bibitem{goulet2017hybridization}
B.~E. Goulet, .F~Roda, and R.~Hopkins.
\newblock Hybridization in plants: old ideas, new techniques.
\newblock {\em Plant Physiology}, 173(1):65--78, 2017.

\bibitem{grant2016introgressive}
P.~R. Grant and B.~R. Grant.
\newblock Introgressive hybridization and natural selection in darwin's
  finches.
\newblock {\em Biological Journal of the Linnean Society}, 117(4):812--822,
  2016.

\bibitem{H73}
J.~A. Hartigan.
\newblock Minimum mutation fits to a given tree.
\newblock {\em Biometrics}, 29:53--65, 1973.

\bibitem{HM06}
K.~T. Huber and V.~Moulton.
\newblock Phylogenetic networks from multi-labelled trees.
\newblock {\em Journal of Mathematical Biology}, 52:613--632, 2006.

\bibitem{huber2019exploring}
K.~T. Huber, V.~Moulton, M.-F. Sagot, and B.~Sinaimeri.
\newblock Exploring and visualizing spaces of tree reconciliations.
\newblock {\em Systematic biology}, 68(4):607--618, 2019.

\bibitem{huber2016folding}
K.~T. Huber, V.~Moulton, M.~Steel, and T.~Wu.
\newblock Folding and unfolding phylogenetic trees and networks.
\newblock {\em Journal of Mathematical Biology}, 73(6-7):1761--1780, 2016.

\bibitem{kundu2018impact}
S.~Kundu and M.~S. Bansal.
\newblock On the impact of uncertain gene tree rooting on
  duplication-transfer-loss reconciliation.
\newblock {\em BMC bioinformatics}, 19(9):21--31, 2018.

\bibitem{mallet2005hybridization}
J.~Mallet.
\newblock Hybridization as an invasion of the genome.
\newblock {\em Trends in Ecology \& Evolution}, 20(5):229--237, 2005.

\bibitem{sankararaman2014genomic}
S.~Sankararaman, S.~Mallick, M.~Dannemann, K.~Pr{\"u}fer, J.~Kelso,
  S.~P{\"a}{\"a}bo, N.~Patterson, and D.~Reich.
\newblock The genomic landscape of neanderthal ancestry in present-day humans.
\newblock {\em Nature}, 507(7492):354--357, 2014.

\bibitem{scholz2019osf}
G.~E. Scholz, A.-A. Popescu, M.~I. Taylor, V.~Moulton, and K.~T. Huber.
\newblock Osf-builder: A new tool for constructing and representing
  evolutionary histories involving introgression.
\newblock {\em Systematic Biology}, 68(5):717--729, 2019.

\bibitem{semple2003phylogenetics}
C.~Semple and M.~Steel.
\newblock {\em Phylogenetics}, volume~24.
\newblock Oxford University Press on Demand, 2003.

\bibitem{solis2016inferring}
C.~Sol{\'\i}s-Lemus and C.~An{\'e}.
\newblock Inferring phylogenetic networks with maximum pseudolikelihood under
  incomplete lineage sorting.
\newblock {\em PLoS Genetics}, 12(3), 2016.

\bibitem{sousa2013understanding}
V.~Sousa and J.~Hey.
\newblock Understanding the origin of species with genome-scale data: modelling
  gene flow.
\newblock {\em Nature Reviews Genetics}, 14(6):404--414, 2013.

\bibitem{stewart2003transgene}
C.~N. Stewart, M.~D. Halfhill, and S.~I. Warwick.
\newblock Transgene introgression from genetically modified crops to their wild
  relatives.
\newblock {\em Nature Reviews Genetics}, 4(10):806--817, 2003.

\bibitem{tofigh2011simultaneous}
A.~Tofigh, M.~Hallett, and J.~Lagergren.
\newblock Simultaneous identification of duplications and lateral gene
  transfers.
\newblock {\em IEEE/ACM Transactions on Computational Biology and
  Bioinformatics}, 8(2):517--535, 2011.

\bibitem{wen2016reticulate}
D.~Wen, Y.~Yu, M.~W. Hahn, and L.~Nakhleh.
\newblock Reticulate evolutionary history and extensive introgression in
  mosquito species revealed by phylogenetic network analysis.
\newblock {\em Molecular Ecology}, 25(11):2361--2372, 2016.

\bibitem{zhang2016tree}
L.~Zhang.
\newblock On tree-based phylogenetic networks.
\newblock {\em Journal of Computational Biology}, 23(7):553--565, 2016.

\bibitem{zhang2016genome}
W.~Zhang, K.~K. Dasmahapatra, J.~Mallet, G.~R.P. Moreira, and M.~R. Kronforst.
\newblock Genome-wide introgression among distantly related heliconius
  butterfly species.
\newblock {\em Genome Biology}, 17(1):25, 2016.

\end{thebibliography}

\end{document}